\documentclass[12pt]{article}
\usepackage{amssymb, amsmath, latexsym, a4}
\usepackage[latin1]{inputenc}
\usepackage[all]{xy}

\begin{document}

\newcommand{\rdg}{\hfill $\Box $}

\newtheorem{De}{Definition}[section]
\newtheorem{Th}[De]{Theorem}
\newtheorem{Pro}[De]{Proposition}
\newtheorem{Le}[De]{Lemma}
\newtheorem{Co}[De]{Corollary}
\newtheorem{Rem}[De]{Remark}
\newtheorem{Ex}[De]{Example}
\newtheorem{Exo}[De]{Exercises}
\newcommand{\tp}{\otimes}
\newcommand{\N}{\mathbb{N}}
\newcommand{\Z}{\mathbb{Z}}
\newcommand{\K}{\mathbb{K}}
\newcommand{\op}{\oplus}
\newcommand{\n}{\underline n}
\newcommand{\es}{{\frak S}}
\newcommand{\ef}{\frak F}
\newcommand{\qu}{\frak Q} \newcommand{\ga}{\frak g}
\newcommand{\la}{\lambda}
\newcommand{\ig}{\frak Y}
\newcommand{\te}{\frak T}
\newcommand{\cok}{{\sf Coker}}
\newcommand{\Hom}{{\sf Hom}}
\newcommand{\im}{{\sf Im}}
\newcommand{\ext}{{\sf Ext}}
\newcommand{\ho}{{\sf H_{AWB}}}
\newcommand{\HH}{{\sf Hoch}}
\newcommand{\adu}{{\rm AWB}^!}

\newcommand{\ele}{\cal L} \newcommand{\as}{\cal A} \newcommand{\ka}{\cal
K}\newcommand{\eme}{\cal M} \newcommand{\pe}{\cal P}

\newcommand{\pn}{\par \noindent}
\newcommand{\pbn}{\par \bigskip \noindent}
\bigskip\bigskip

\centerline {\Large {\bf On universal central extensions of Hom-Lie
algebras}}

\

\centerline {\bf J. M. Casas${^1}$, M. A. Insua${^1}$ and N. Pacheco${^2}$}

\bigskip \bigskip
\centerline{${^1}$ Dpto.  Matem\'atica Aplicada I, Univ. de Vigo, 36005 Pontevedra, Spain}
\centerline{e-mail addresses: jmcasas@uvigo.es, avelino.ainsua@gmail.com}

 \bigskip
 \centerline{$^{2}$ IPCA, Dpto. de Ciências, Campus do IPCA,
 Lugar do Aldão}
\centerline{4750-810 Vila Frescainha, S. Martinho, Barcelos,
 Portugal}
\centerline{e-mail address: natarego@gmail.com}
\bigskip \bigskip \bigskip \bigskip

\par
{\bf Abstract:} We develop a theory of universal central extensions of Hom-Lie algebras. Classical results of universal central extensions of Lie algebras cannot be completely extended to Hom-Lie algebras setting,  because of the composition of two central extensions is not central. This fact leads to introduce the notion of universal $\alpha$-central extension. Classical results as the existence of a universal central extension of a perfect Hom-Lie algebra remains true, but others as the central extensions of the middle term of a universal central extension  is split only holds for $\alpha$-central extensions.
A homological characterization of universal ($\alpha$)-central extensions is given.

\bigskip \bigskip


\bigskip \bigskip
{\it A. M. S. Subject Class. (2010):} 17A30, 17B55, 17B60, 18G35, 18G60

\section{Introduction}
The Hom-Lie algebra structure was initially introduced in \cite{HLS}  motivated by examples of deformed Lie algebras coming from twisted discretizations of vector fields. Hom-Lie algebras are $\mathbb{K}$-vector spaces endowed with a bilinear skew-symmetric bracket satisfying a  Jacobi identity twisted by a map. When this map is the identity map, then the definition of Lie algebra is recovered.

The study of this algebraic structure was the subject of several papers \cite{HLS, MS, MS2, Sh, Yau}. In particular, a homology theory for Hom-Lie algebras, which generalizes the Chevalley-Eilenberg homology for a Lie algebra, was the subject of \cite{Yau1, Yau2}.

In the classical setting, homology theory is closely related with universal central extensions. Namely, the second homology with trivial coefficients group is the kernel of the universal central extension and universal central extensions are characterized by means of the first and second homologies with trivial coefficients.

Our goal in the present paper is to investigate if the homology for Hom-Lie algebras introduced in \cite{Yau1, Yau2} allows the characterization of universal central extensions of Hom-Lie algebras in terms of Hom-Lie homologies.
But when we try to generalize the classical results of universal  central extensions theory
of Lie algebras to Hom-Lie algebras  an important problem occurs, namely the composition of central extensions
is not central in general. This fact  doesn't allow  a complete generalization
of classical results, however requires the introduction of a new concept of centrality for Hom-Lie algebra extensions.

To show our results, we organize the paper as follows: in Section 2 we recall some basic needed material on Hom-Lie algebras, the notions of center, commutator and module. In order to have examples, we include the classification of two-dimensional complex Hom-Lie algebras. In section 3 we recall the chain complex given in \cite{Yau2} and we prove its well-definition by means of the Generalized Cartan's formulas; the interpretation of low-dimensional homologies is given. In section 4 we present our main results on universal central extensions, namely we extend classical results and present a counterexample showing that the composition of two central extension is not a central extension. This fact lead us to define $\alpha$-central  extensions as extensions for which the image by the twisting endomorphism $\alpha$ of the kernel is included in the center of the middle Hom-Lie algebra. We can extend classical results as: a Hom-Lie algebra is perfect if and only if admits a universal central extension and the kernel of the universal central extension is the second homology with trivial coefficients of the Hom-Lie algebra. Nevertheless, other result as: if a central extension $0 \to (M, \alpha_M) \stackrel{i} \to (K,\alpha_K) \stackrel{\pi} \to (L, \alpha_L) \to 0$ is universal, then $(K,\alpha_K)$ is perfect and every central extension of $(K,\alpha_K)$ is split only holds for universal $\alpha$-central extensions, which means that only lifts on $\alpha$-central extensions. Other relevant result, which cannot be extended in the usual way, is: if $0 \to (M, \alpha_M) \stackrel{i} \to (K,\alpha_K) \stackrel{\pi} \to (L, \alpha_L) \to 0$ is a universal   $\alpha$-central extension, then $H_1^{\alpha}(K) = H_2^{\alpha}(K) = 0$. Of course, when the twisting endomorphism is the identity morphism, then all the new notions and all the new results coincide with the classical ones.

\section{Hom-Lie algebras}

\begin{De} \cite{HLS} \label{def}
A Hom-Lie algebra is a triple $(L,[-,-],\alpha_L)$ consisiting of a $\mathbb{K}$-vector space $L$, a bilinear map $[-,-] : L \times L \to L$ and a $\mathbb{K}$-linear map $\alpha_L : L \to L$ satisfying:
\begin{enumerate}
\item [a)] $[x,y] = - [y,x]$ \ \ \ \ \ \ \ \ \ \ \ \ \ \ \ \ \ \ \ \ \ \ \ \ \ \ \ \ \ \ \ \ \ \ \ \ \ \ \ \ \ \ \ (skew-symmetry)
\item [b)] $[\alpha_L(x),[y,z]]+[\alpha_L(z),[x,y]]+[\alpha_L(y),[z,x]]=0$ \ \ (Hom-Jacobi identity)
\end{enumerate}
for all $x, y, z \in L$.
\end{De}

In terms of the adjoint representation $ad_x : L  \to L, ad_x(y)=[x,y]$, the Hom-Jacobi identity can be written as follows \cite{MS}: $$ad_{{\alpha_L}(z)} \ ad_y = ad_{{\alpha_L}(y)} \ ad_z + ad_{[z,y]} \ \alpha_L$$

\begin{De} \cite{Yau1}
A Hom-Lie algebra $(L,[-,-],\alpha_L)$ is said to be multiplicative  if the linear map $\alpha_L$ preserves the bracket.
\end{De}

\begin{Ex}\label{ejemplo 1} \
\begin{enumerate}
\item[a)] Taking $\alpha_L = Id$ in Definition \ref{def} we obtain the definition of a Lie algebra. Hence Hom-Lie algebras include Lie algebras as a subcategory, thereby motivating the name "Hom-Lie algebras" as a deformation of Lie algebras twisted by an endomorphism. Moreover it  is a multiplicative Hom-Lie algebra.

    \item[b)] Let $(A,\mu_A, \alpha_A)$ be a multiplicative Hom-associative algebra \cite{MS}. Then $HLie(A)=(A,[-,-],\alpha_A)$ is a multiplicative Hom-Lie algebra in which $[x,y]=\mu_A(x,y)-\mu_A(y,x)$, for all $x,y \in A$ \cite{Yau1}.

  \item[c)] Let $(L,[-,-])$ be a Lie algebra and $\alpha:L \to L$ be a Lie algebra endomorphism. Define $[-,-]_{\alpha} : L \otimes L \to L$ by $[x,y]_{\alpha} = \alpha[x,y]$, for all $x, y \in L$. Then $(L,[-,-]_{\alpha}, \alpha)$ is a multiplicative Hom-Lie algebra \cite[Th. 5.3]{Yau1}.

\item[d)] Abelian or commutative Hom-Lie algebras are $\mathbb{K}$-vector spaces $V$ with trivial bracket and any linear map $\alpha:V \to V$ \cite{HLS}.

\item[e)] The Jackson Hom-Lie algebra  $\frak{sl}_2(\mathbb{K})$ is a Hom-Lie deformation of the classical Lie algebra  $\frak{sl}_2(\mathbb{K})$ defined by $[h,f]=-2f, [h,e]=2e,[e,f]=h$. The Jackson $\frak{sl}_2(\mathbb{K})$ is related to derivations. As a $\mathbb{K}$-vector space is generated by $e, f, h$ with multiplication given by $[h,j]_t=-2f-2tf, [h,e]_t=2e, [e,f]_t=h+ \frac{t}{2} h$ and the linear map $\alpha_t$ is defined by $\alpha_t(e) = \frac{2+t}{2(1+t)}e= e+ \displaystyle \sum_{k=0}^{\infty} \frac{(-1)^k}{2} t^k e, \alpha_t(h)=h, \alpha_t(f)=f+ \frac{t}{2} f$ \cite{MS2}.

    \item[f)] For examples coming from deformations we refer to \cite{Yau1}.
\end{enumerate}
\end{Ex}

\begin{De}\label{homo}
A homomorphism of  Hom-Lie algebras $f:(L,[-,-],\alpha_L) \to (L',[-,-]',\alpha_{L'})$ is a $\mathbb{K}$-linear map $f : L \to L'$ such that
\begin{enumerate}
\item[a)] $f([x,y]) =[f(x),f(y)]'$
\item [b)] $f \ \alpha_L(x) = \alpha_{L'} \ f(x)$
\end{enumerate}
for all $x, y \in L$.

The Hom-Lie algebras $L$ and $L'$ are isomorphic if there is a bijective Hom-Lie algebras homomorphism $f: L \to L'$.
\end{De}

A homomorphism of multiplicative Hom-Lie algebras is a homomorphism of the underlying Hom-Lie algebras.
\bigskip

So we have defined the category ${\sf Hom-Lie}$ (respectively, $({\sf Hom-Lie_{\rm mult}})$ whose objects are Hom-Lie (respectively, multiplicative Hom-Lie) algebras and whose morphisms are the homomorphisms of Hom-Lie (respectively, multiplicative Hom-Lie) algebras.
There is an obvious inclusion functor $inc : {\sf Hom-Lie_{\rm mult}} \to {\sf Hom-Lie}$. This functor has as left adjoint the multiplicative functor $(-)_{\rm mult} : {\sf Hom-Lie} \to {\sf Hom-Lie_{\rm mult}}$ which assigns to a Hom-Lie algebra $(L,[-,-],\alpha_L)$ the Hom-Lie multiplicative algebra $(L/I,[-,-],\overline{\alpha})$, where $I$ is the ideal of $L$ generated by the elements $\alpha_L[x,y]-[\alpha_L(x),\alpha_L(y)]$, for all $x, y \in L$ and $\overline{\alpha}$ is induced by $\alpha$.
\bigskip

In the sequel we refer Hom-Lie algebra to a multiplicative Hom-Lie algebra.
\bigskip

Let $(L, [-,-],\alpha_L)$  be an $n$-dimensional Hom-Lie algebra with basis $\{a_1, a_2, \dots,a_n \}$ and endomorphism $\alpha_L$ represented by the matrix $A = (\alpha_{ij})$ with respect to the given basis. To determine its algebraic structure is enough to know its structural constants, i.e. the scalars $c_{ij}^k$   such that $[a_{i}, a_{j}] = \displaystyle \sum_{k=1}^n c_{ij}^k a_k$, and the entries  $\alpha_{ij}$ corresponding to the matrix $A$. These terms are related according to the following

\begin{Pro} \label{constantes estructura HomLie}  Let  $(L, [-,-],\alpha_L)$ be a Hom-Lie algebra with basis $\{a_1, a_2, \dots,a_n \}$.  Let $c_{ij}^k, 1 \leq i,j,k \leq n$ be the  structural constants relative to this basis and $\alpha_{ij}, 1 \leq i, j \leq n$ the entries of the  matrix $A$ associated to the endomorphism $\alpha_L$ with  respect to the given basis. Then $(L, [-,-],\alpha_L)$  is a Hom-Lie algebra if and only if the structural constants and the  matricial entries satisfy the following properties:
\begin{enumerate}
 \item[a)] $c_{ij}^k+c_{ji}^k=0, 1 \leq i,j,k \leq n ;\quad  c_{ii}^k =0, 1\leq i, k \leq n$, char$(\mathbb{K}) \neq 2$.
  \item[b)]  $ \displaystyle \sum_{p=1}^n \alpha_{pi} \left( \sum_{q=1}^n c_{jk}^q c_{pq}^l \right) + \sum_{p=1}^n \alpha_{pk} \left( \sum_{q=1}^n c_{ij}^q c_{pq}^l \right) + \sum_{p=1}^n \alpha_{pj} \left( \sum_{q=1}^n c_{ki}^q c_{pq}^l \right)  = 0,$

$1 \leq i, j, k, l, \leq n.$
  \end{enumerate}
\end{Pro}
{\it Proof.}\

\noindent a) There is not difference with Lie-algebras case \cite{Hu}.
\medskip

\noindent b) Applying Hom-Jacobi identity \ref{def} b):
$$[\alpha(a_i),[a_j,a_k]] + [\alpha(a_k),[a_i,a_j]] + [\alpha(a_j),[a_k,a_i]] = 0$$
$$[\sum_{l=1}^n \alpha_{li} a_l, \sum_{m=1}^n c_{jk}^m a_m]+[\sum_{p=1}^n \alpha_{pk} a_p, \sum_{q=1}^n c_{ij}^q a_q]+ [\sum_{r=1}^n \alpha_{rj} a_r, \sum_{s=1}^n c_{ki}^s a_s] = 0$$
$$\sum_{p=1}^n \alpha_{pi} \left( \sum_{q=1}^n c_{jk}^q [a_p,a_q] \right) + \sum_{p=1}^n \alpha_{pk} \left( \sum_{q=1}^n c_{ij}^q [a_p,a_q] \right) + \sum_{p=1}^n \alpha_{pj} \left( \sum_{q=1}^n c_{ki}^q [a_p,a_q] \right) =0 $$
$$\sum_{l=1}^n \left\{ \sum_{p=1}^n \alpha_{pi} \left( \sum_{q=1}^n c_{jk}^q c_{pq}^l \right) + \sum_{p=1}^n \alpha_{pk} \left( \sum_{q=1}^n c_{ij}^q c_{pq}^l \right) + \sum_{p=1}^n \alpha_{pj} \left( \sum_{q=1}^n c_{ki}^q c_{pq}^l \right) \right\} a_l = 0$$

 \rdg

\begin{Pro} \label{isom}
The Hom-Lie algebras $(L,[-,-],\alpha_L)$ and $(L,[-,-]',\alpha_{L'})$ with same underlying $\mathbb{K}$-vector space are isomorphic if and only if there exists a regular matrix $P$ such that $A' = P^{-1}. A. P$ and $P.[a_i,a_j] = [P.a_i, P.a_j]'$, where $A, A'$ and $P$ denote the corresponding matrices representing $\alpha_L, \alpha_{L'}$ and $f$ with respect to the basis $\{a_1, \dots, a_n\}$, respectively.
\end{Pro}
{\it Proof.} The fact comes directly from Definition \ref{homo}. \rdg

\begin{Pro}\label{clas HomLie}
The 2-dimensional complex multiplicative Hom-Lie algebras with basis $\{a_1, a_2\}$ are isomorphic to one in the following isomorphism classes:
\begin{enumerate}
\item[a)] Abelian.
\item[b)] $[a_1,a_2] = -[a_2,a_1] = a_1$ and $\alpha$ is represented by the matrix $\left( \begin{array}{cc} 0 & \alpha_{12} \\ 0 & \alpha_{22} \end{array} \right)$.
\item[c)] $[a_1,a_2] = -[a_2,a_1] = a_1$ and $\alpha$ is represented by the matrix $\left( \begin{array}{cc} \alpha_{11} & \alpha_{12} \\ 0 & 1 \end{array} \right)$, with $\alpha_{11} \neq 0$.
\end{enumerate}
\end{Pro}
{\it Proof.} From the skew-symmetry condition we have that $[a_1,a_1]=[a_2,a_2]=0$ and $[a_1,a_2]=-[a_2,a_1]$. The Hom-Jacobi identity \ref{def} $b)$ is satisfied independently of the homomorphism $\alpha$. So we only have restrictions coming from the fact that the $\mathbb{C}$-linear map $\alpha_L : L \to L$ represented by the matrix $\left( \begin{array}{cc} \alpha_{11} & \alpha_{12} \\ \alpha_{21} & \alpha_{22} \end{array} \right)$ must preserve the bracket.

\noindent First at all, we apply the change of basis given by the equations $\left\{ \begin{array}{ll} a_1' =x.a_1 + y.a_2 \\ a_2'= \frac{1}{x}.a_2 \end{array} \right.$, if $x \neq 0$, and $\left\{ \begin{array}{ll} a_1' = a_2 \\ a_2'= -\frac{1}{y}.a_1 \end{array} \right.$, if $x=0$ and $y \neq 0$, to normalize the bracket, obtaining the bracket $[a_1',a_1']=[a_2',a_2']=0, [a_1',a_2']=-[a_2',a_1']=p.a_1'$, for $p=0, 1$.

\noindent From the fact that $\alpha_L : L \to L$ preserves the bracket, we derive the following equations:
$$\left. \begin{array}{rcl} (\alpha_{11} \alpha_{22} - \alpha_{21} \alpha_{12}). p &= & p. \alpha_{11} \\ p. \alpha_{21} & = & 0 \end{array} \right\}$$
which reduces to the following system:
$$\left. \begin{array}{rcl} p. \alpha_{11}.( \alpha_{22} - 1) &= & 0 \\ p. \alpha_{21} & = & 0 \end{array} \right\}$$
Hence, for $p=0$ the system is trivially satisfied. All the matrices representing $\alpha$ are valid and the bracket is trivial, so $(L, [-,-], \alpha_L)$ is an abelian Hom-Lie algebra.
In case $p=1$, we derive the matrices corresponding to the cases {\it b)} and {\it c)}.

The different classes obtained are not pairwise isomorphic thanks to Proposition \ref{isom}.
\rdg

\begin{Rem}\
\begin{enumerate}

\item[a)] Two algebras of the class  $b)$ in Proposición \ref{clas HomLie}, with endomorphisms given by the matrices   $\left( \begin{array}{cc} 0 & \alpha_{12} \\ 0 & \alpha_{22} \end{array} \right)$ and $\left( \begin{array}{cc} 0 & \beta_{12} \\ 0 & \beta_{22} \end{array} \right)$, are isomorphic if and only if $\alpha_{22}=\beta_{22}$ and $\beta_{12}= p. \alpha_{12} + q. \alpha_{22}, p, q \in \mathbb{C}, p \neq 0$.

\item[b)] Two algebras of the class  $c)$ in Proposición \ref{clas HomLie}, with endomorphisms given by the matrices   $\left( \begin{array}{cc} \alpha_{11} & \alpha_{12} \\ 0 & 1 \end{array} \right)$ and $\left( \begin{array}{cc} \beta_{11} & \beta_{12} \\ 0 & 1 \end{array} \right)$, are isomorphic if and only if $\alpha_{11}=\beta_{11}$ and $\beta_{12}= p. \alpha_{12} - q. \alpha_{11}+q, p, q \in \mathbb{C}, p \neq 0$.

\item[c)] Obviously if $\Phi : (L,[-,-], \alpha_L) \to (L,[-,-]',\alpha_{L'})$ is an isomorphism of Hom-Lie algebras, then $det(\alpha_L) = det(\alpha'_L)$. Consequently, if $det(\alpha_L) \neq det(\alpha'_L)$, then the Hom-Lie algebras are not isomorphic.
\item[d)] The following  table shows  by means of its algebraic properties that the classes given in Proposición \ref{clas HomLie} are not pairwise isomorphic.

\begin{center}

\begin{tabular}{|l|c|c|}
\hline   & Abelian & $det(\alpha)$ \\
\hline \ref{clas HomLie} a) & Yes &  \\
\hline \ref{clas HomLie} b) & Non & 0 \\
 \hline \ref{clas HomLie} c) & Non & $\neq 0$  \\

\hline
\end{tabular}

Complex two-dimensional Hom-Lie algebras
\end{center}
\end{enumerate}
\end{Rem}

\begin{De}
Let $(L,[-,-],\alpha_L)$ be a Hom-Lie algebra. A  Hom-Lie subalgebra $H$ is a linear subspace of $L$, which is closed for the bracket and invariant by $\alpha$, that is,
\begin{enumerate}
\item [a)] $[x,y] \in H,$ for all $x, y \in H$
\item [b)] $\alpha(x) \in H$, for all $x \in H$
\end{enumerate}

A  Hom-Lie subalgebra $H$ of $L$ is said to be a  Hom-ideal if $[x, y] \in H$ for all $x \in H, y \in L$.

If $H$ is a  Hom-ideal of $L$, then  $(L/H, [-,-], \overline{\alpha_L})$ naturally inherits a structure of Hom-Lie algebra, which is said to be the quotient Hom-Lie algebra.
\end{De}

\begin{De}
Let $H$ and $K$ be  Hom-ideals of a Hom-Lie algebra $(L,[-,-],\alpha_L)$. The commutator  Hom-Lie subalgebra of $H$ and $K$, denoted by $[H,K]$, is the  Hom-subalgebra of $L$ spanned by the brackets $[h,k], h \in H, k \in K$.
\end{De}

\begin{Le}\label{ideales} Let $H$ and $K$ be  Hom-ideals of a Hom-Lie algebra
$(L,[-,-],\alpha_L)$. The following statement hold:

\begin{enumerate} \label{ideales}

\item[a)] $H \cap K$ and $H+K$ are  Hom-ideals of $L$.

\item[b)]  $[H,K] \subseteq H \cap K$.

\item[c)] $[H,K]$ is a Hom-ideal of $L$ when $\alpha_L$ is surjective.

\item[d)] $[H,K]$ is a Hom-ideal of $H$ and $K$, respectively.

\item[f)]  If $H=K=L$, then $[L,L]$ is a Hom-ideal of $L$.
\end{enumerate}
\end{Le}

\begin{Le}
Let $H$ and $K$ be Hom-ideals of a Hom-Lie algebra
$(L,[-,-],\alpha_L)$, then $[H,K]$ is a Hom-ideal of $(\alpha_L(L),[-,-],\alpha_{L_{\mid}})$.
\end{Le}

\begin{De}
The center of a  Hom-Lie algebra  $(L,[-,-],\alpha_L)$ is the $\mathbb{K}$-vector subspace  $$Z(L) = \{ x \in L \mid [x, y] =0 , {\rm for\ all}\ y \in L\}$$
\end{De}

\begin{Rem}
When $\alpha_L : L \to L$ is a surjective endomorphism, then  $Z(L)$ is a Hom-ideal of $L$.
\end{Rem}

 \begin{De} \label{action}
 Let $(L,[-,-],\alpha_L)$ and $(M,[-,-], \alpha_M)$ be  Hom-Lie algebras. A Hom-$L$-action from $(L,[-,-],\alpha_L)$ over  $(M,[-,-], \alpha_M)$ consists in  a bilinear map $\rho: L \otimes M \to M,$ given by $ \rho(x \otimes m)=x \centerdot m,$ satisfying the following properties:
 \begin{enumerate}
 \item[a)] $[x,y] \centerdot \alpha_M(m) = \alpha_L(x) \centerdot (y \centerdot m) - \alpha_L(y) \centerdot (x \centerdot m)$
 \item [b)] $\alpha_L(x) \centerdot[m,m'] = [x \centerdot m, \alpha_M(m')]+[\alpha_M (m), x \centerdot m']$
 \item [c)] $\alpha_M(x \centerdot m) = \alpha_L(x) \centerdot \alpha_M(m)$
 \end{enumerate}
 for all $x, y \in L$ and $m, m' \in M$.

 Under these circumstances, we say that $(L, \alpha_L)$ Hom-acts over $(M, \alpha_M)$.
 \end{De}

 \begin{Rem}
 When $(M,[-,-], \alpha_M)$ is an abelian Hom-Lie algebra, Definition \ref{action} goes back to the definiton of Hom-L-module in  \cite{Yau1}.
 \end{Rem}

\begin{Ex}\label{ejemplo 2} \
\begin{enumerate}
\item[a)] $(L,[-,-],\alpha_L)$ acts on itself by the action given by the bracket.
\item[b)] Let  $\frak{g}$ and $\frak{m}$ be   Lie algebras with a Lie action from $\frak{g}$ over  $\frak{m}$. Then $(\frak{g}, Id_{\frak{g}})$ Hom-acts over $(\frak{m}, Id_{\frak{m}})$.
    \item[c)] Let $\frak{g}$ be a Lie algebra, $\alpha : \frak{g} \to \frak{g}$ an endomorphism and $M$ an $\frak{g}$-module in the usual sense, such that the action from $\frak{g}$ over $M$ satisfies the condition $\alpha(x) \centerdot m = x \centerdot m$, for all $x \in \frak{g}$ and $m \in M$. Then $(M, Id)$ is a Hom-$\frak{g}$-module.

        An example of this situation is given by the 2-dimensional Lie algebra $L$ generated by $\{ e,f \}$ with bracket $[e,f]=-[f,e]=e$ and endomorphism $\alpha$ represented by the matrix $\left( \begin{array}{cc} 1 & 1 \\ 0 & 1 \end{array} \right)$, where $M$ the ideal generated by $\{e\}$.
  \item [d)] An abelian sequence of Hom-Lie algebras is an exact sequence of Hom-Lie algebras $0 \to (M,\alpha_M) \stackrel{i}\to (K,\alpha_K) \stackrel{\pi}\to (L,\alpha_L) \to 0$, where $(M,\alpha_M)$ is an  abelian Hom-Lie algebra, $\alpha_K \ i = i  \ \alpha_M$ and $\pi \ \alpha_K = \alpha_L \ \pi$.

      The abelian sequence induces a Hom-$L$-module structure on $(M, \alpha_M)$ by means of the action given by $\rho : L \otimes M \to M, \rho(l,m)=[k,m], \pi(k)=l$.

\item[e)] For other examples we refer to Example 6.2 in \cite{Yau1}.
\end{enumerate}
\end{Ex}

\section{Homology}

Following \cite{Yau1, Yau2}, for a Hom-Lie algebra $(L, [-,-], \alpha_{L})$ and a (right) Hom-$L$-module  $\left(  M,\alpha_{M}\right)$, one denotes by
$$C_{n}^{\alpha}\left(  L,M\right)  :=M\otimes \Lambda^n L,\ n\geqslant0$$
the $n$-chain module of $(L, [-,-], \alpha_{L})$ with coefficients in $\left(  M,\alpha_{M}\right)$.

For $n\geqslant1$, one defines the   $\mathbb{K}$-linear map,
$$d_{n}:C_{n}^{\alpha}\left(  L,M\right)  \longrightarrow C_{n-1}^{\alpha}\left(  L,M\right)$$
by
$$d_{n}\left(  m\otimes x_{1}\wedge\cdots\wedge
x_{n}\right)  =\overset{n}{\underset{i=1}{%
{\displaystyle\sum}
}}\left(  -1\right)  ^{i+1}m\centerdot x_{i}\otimes\alpha_{L}(x_{1})
\wedge\cdots\wedge\widehat{\alpha_{L}(x_{i})}\wedge\cdots\wedge \alpha_{L}(x_{n}) +$$
$$\underset{1\leqslant i<j\leqslant n}{\sum}\left(  -1\right)
^{i+j}\alpha_{M}\left(  m\right)  \otimes\left[  x_{i},x_{j}\right]
\wedge\alpha_{L}\left(  x_{1}\right)  \wedge\cdots\wedge\widehat{\alpha
_{L}\left(  x_{i}\right)  }\wedge\cdots\wedge\widehat{\alpha_{L}\left(
x_{j}\right)  }\wedge\cdots\wedge\alpha_{L}\left(  x_{n}\right)$$

 Although in \cite{Yau1, Yau2} is proved  that  $(CL_{n}^{\alpha}\left(  L,M\right),d_n)$ is a well-defined chain complex, we present an alternative  proof by means of a generalization of Cartan's formulas.
Firstly, we define for all $y\in L$ and $n\in \mathbb{N}$,  two linear maps,
$$\theta_{n}\left(  y\right)
:C_{n}^{\alpha}\left(  L,M\right)  \longrightarrow C_{n}^{\alpha}\left(
L,M\right)$$ by
$$\theta_{n}\left(  y\right)  \left(  m\otimes
x_{1}\wedge\cdots\wedge x_{n}\right)  =-m\centerdot y\otimes\alpha_{L}\left(
x_{1}\right)  \wedge\cdots\wedge\alpha_{L}\left(  x_{n}\right) +$$
$$\overset{n}{\underset{i=1}{%
{\displaystyle\sum}
}}\left(  -1\right)  ^{i}\alpha_{M}\left(  m\right)  \otimes\left[
x_{i},y\right]  \wedge\alpha_{L}\left(  x_{1}\right)  \wedge\cdots
\wedge\widehat{\alpha_{L}\left(  x_{i}\right)  }\wedge\cdots\wedge\alpha
_{L}\left(  x_{n}\right)$$
and
$$i_{n}\left( \alpha_L(y) \right)
:C_{n}^{\alpha}\left(  L,M\right)  \longrightarrow C_{n+1}^{\alpha}\left(
L,M\right)$$ by
$$i_{n}\left(  \alpha_{L}\left(
y\right)  \right)  \left(  m\otimes x_{1}\wedge\cdots\wedge x_{n}\right)
=\left(  -1\right)  ^{n}m\otimes x_{1}\wedge\cdots\wedge x_{n}\wedge y$$

\begin{Pro} (Generalized Cartan's formulas) \label{Cartan Lie}

The following identities hold:
\begin{enumerate}

\item[a)] $d_{n+1} i_{n}\left(  \alpha_{L}\left(  y\right)  \right)
+i_{n-1}\left(  \alpha_{L}^{2}\left(  y\right)  \right)  d_{n}%
=-\theta_{n}\left(  y\right)  $, for all  $n\geqslant1$.

\item[b)] $\theta_{n}\left(  \alpha_{L}\left(  x\right)  \right)
\theta_{n}\left(  y\right)  -\theta_{n}\left(  \alpha_{L}\left(  y\right)
\right)  \theta_{n}\left(  x\right)  =\theta_{n}\left(  \left[
x,y\right]  \right)  \left(  \alpha_{M}\otimes\alpha_{L}^{\wedge n}\right)  $,
for all $n \geq 0$.

\item[c)] $\theta_{n}\left(  x\right)   i_{n-1}\left(  \alpha
_{L}\left(  y\right)  \right)  -i_{n-1}\left(  \alpha_{L}^{2}\left(  y\right)
\right)  \theta_{n-1}\left(  x\right)  =  i_{n-1}(\alpha_{L}\left[
x,y\right])  \left(  \alpha_{M}\otimes\alpha_{L}^{\wedge (n-1)}\right)  $, for all
$n \geq 1$.

\item[d)] $\theta_{n-1}\left(  \alpha_{L}\left(  y\right)  \right)
d_{n}=d_{n} \theta_{n}\left(  y\right)$, for all
$n \geq 1$.

\item[e)] $d_{n}  d_{n+1}=0,$ for all $n \geq 1$.

\end{enumerate}
\end{Pro}
{\it Proof.} The proof follows with a routine induction, so we omit it. \rdg
\bigskip

In case $\alpha_L = Id_L, \alpha_M = Id_M$, the above formulas become to the Cartan's formulas for the Chevalley-Eilenberg homology \cite{HS}.
\bigskip

Thanks to Proposition \ref{Cartan Lie},  $\left(  C_{\star}^{\alpha}\left(  L,M\right), d_{\star}
\right)$ is a well-defined chain complex (an alternative proof can be seen in  \cite{Yau2}). Its homology is said to be the
homology of the  Hom-Lie algebra $\left(  L,\left[  -,-\right]  ,\alpha
_{L}\right)$ with coefficients in the Hom-$L$-module $\left(  M,\alpha
_{M}\right)$ and it is denoted by:
$$ H_{\star}^{\alpha
}\left(  L,M\right)  :=H_{\star}\left(  C_{\star}^{\alpha}\left(  L,M\right)
,d_{\star}\right)$$

An easy computation in low-dimensional  cycles and boundaries provides the following results:
$$H_{0}^{\alpha
}\left(  L,M\right)  =\frac{Ker\left(  d_{0}\right)  }{\operatorname{Im}%
\left(  d_{1}\right)  }=\frac{M}{M_{L}}$$
where $M_{L}=\left\{  m\centerdot l:m\in M,l\in L\right\}$.

Now let us consider $M$ as a trivial Hom-$L$-module, i.e.  $m \centerdot l=0$, then
$$H_{1}^{\alpha
}\left(  L,M\right)=\frac{Ker\left(  d_{1}\right)  }{\operatorname{Im}%
\left(  d_{2}\right)  }  =\frac{M\otimes L}{\alpha\left(  M\right)  \otimes\left[
L,L\right]  }$$
In particular, if $M = \mathbb{K}$, then $H_{1}^{\alpha
}\left(  L,\mathbb{K} \right)  =\frac{L}{\left[
L,L\right]  }$.
\bigskip

\section{Universal central extensions}

Through this section we will deal with universal central extensions
of  Hom-Lie algebras. We will generalize classical results of universal  central extensions theory
of Lie algebras, but here an important problem appears, namely the composition of central extensions
is not central in general,  as the Example \ref{contraejemplo} shows. This fact doesn't allow  a complete generalization
of classical results, however requires the introduction of a new concept of centrality for Hom-Lie algebra extensions.

\begin{De} \label{alfacentral}
 A short exact sequence of Hom-Lie algebras  $(K) : 0 \to (M, \alpha_M) \stackrel{i} \to (K,\alpha_K) \stackrel{\pi} \to (L, \alpha_L) \to 0$ is said to be central if $[M, K] = 0 $. Equivalently,  $M \subseteq Z(K)$.

The sequence $(K)$ is said to be $\alpha$-central if $[\alpha_M(M), K] = 0$ . Equivalently, $\alpha_M(M) \subseteq Z(K)$.
\end{De}

\begin{Rem}
Let us observe that  both notions coincide when $\alpha_M = Id_M$. Obviously, every  central extension is an $\alpha$-central extension, but the converse doesn't hold as the following counterexample shows:

Consider the two-dimensional Hom-Lie algebra $L$ with basis $\{ a_1, a_2 \}$, bracket given by $$[a_1, a_2] = - [a_2, a_1] = a_1$$
and endomorphism $\alpha_L = 0$.

Let be the three-dimensional Hom-Lie algebra $K$ with basis $\{b_1, b_2, b_3 \}$, bracket given by
$$[b_1,b_2] = - [b_2, b_1] = b_1; [b_1,b_3] = - [b_3, b_1] = b_1; [b_2,b_3] = - [b_3, b_2] = b_2$$
and endomorphism $\alpha_K = 0$.

The surjective homomorphism $\pi : (K, 0) \to (L, 0)$ given by $\pi(b_1) = 0, \pi(b_2) = a_1, \pi(b_3) = a_2$ is an $\alpha$-central extension, since Ker $(\pi) = \langle \{b_1 \} \rangle$ and $[\alpha_K(Ker (\pi)), K] = 0$, but is not a central extension, since $[Ker (\pi), K] = \langle \{ b_1 \} \rangle$.
\end{Rem}

\begin{De}
A central extension  $(K) : 0 \to (M, \alpha_M) \stackrel{i} \to (K,\alpha_K) \stackrel{\pi} \to (L, \alpha_L) \to 0$ is said to be universal if for every   central extension $(K) : 0 \to (M', \alpha_{M'}) \stackrel{i'} \to (K',\alpha_{K'}) \stackrel{\pi'} \to (L, \alpha_L) \to 0$  there exists a  unique homomorphism of Hom-Lie algebras  $h : (K,\alpha_K) \to (K',\alpha_{K'})$ such that $\pi'\ h = \pi$.

 A central extension  $(K) : 0 \to (M, \alpha_M) \stackrel{i} \to (K,\alpha_K) \stackrel{\pi} \to (L, \alpha_L) \to 0$ is said to be universal $\alpha$-central  if for every  $\alpha$-central extension $(K) : 0 \to (M', \alpha_{M'}) \stackrel{i'} \to (K',\alpha_{K'}) \stackrel{\pi'} \to (L, \alpha_L) \to 0$  there exists a unique homomorphism of  Hom-Lie algebras  $h : (K,\alpha_K) \to (K',\alpha_{K'})$ such that $\pi' \ h = \pi$.
\end{De}

\begin{Rem} \label{rem}
Obviously, every universal $\alpha$-central extension  is a  universal  central extension.
Let us observe that both notions coincide when $\alpha_M = Id_M$.
\end{Rem}

\begin{De}
A Hom-Lie algebra  $(L, \alpha_L)$ is said to be perfect if  $L = [L, L]$.
\end{De}

\begin{Le} \label{lema 1}
Let $\pi : (K,\alpha_K) \to (L, \alpha_L)$ be a surjective homomorphism  of  Hom-Lie algebras. If $(K,\alpha_K)$ is a perfect  Hom-Lie algebra, then $(L, \alpha_L)$ also it is.
\end{Le}

\begin{Le} \label{lema 2}
Let $0 \to (M, \alpha_M) \stackrel{i} \to (K,\alpha_K) \stackrel{\pi} \to (L, \alpha_L) \to 0$ be a central extension and $(K,\alpha_K)$ a perfect Hom-Lie algebra. If there exists a homomorphism of Hom-Lie algebras  $f : (K,\alpha_K) \to (A, \alpha_A)$ such that $\tau \ f = \pi$, where $0 \to (N, \alpha_N) \stackrel{j} \to (A,\alpha_A) \stackrel{\tau} \to (L, \alpha_L) \to 0$ is a central extension, then $f$ is unique.
\end{Le}

The proofs of these two last Lemmas use classical arguments, so we omit it.

\begin{Le} \label{lema 3}
If $0 \to (M, \alpha_M) \stackrel{i} \to (K,\alpha_K) \stackrel{\pi} \to (L, \alpha_L) \to 0$ is a universal central extension, then  $(K,\alpha_K)$ and $(L, \alpha_L)$ are perfect Hom-Lie algebras.
\end{Le}
{\it Proof.} Let us assume that $(K,\alpha_K)$ is not a perfect Hom-Lie, then $[K,K] \varsubsetneq K$. Hence $(K/[K,K], \widetilde{\alpha})$, where $\widetilde{\alpha}$ is the induced  homomorphism, is an abelian Hom-Lie algebra, consequently, it  is a trivial Hom-L-module. Let us consider the central extension $0 \to (K/[K,K], \widetilde{\alpha}) \to (K/[K,K] \times L, \widetilde{\alpha}\times \alpha_L) \stackrel{pr}\to  (L, \alpha_L) \to 0$. Then the  homomorphisms of Hom-Lie algebras  $\varphi, \psi : (K,\alpha_K) \to (K/[K,K] \times L, \widetilde{\alpha}\times \alpha_L)$ given by $\varphi(k)=(k+[K,K],\pi(k))$ and $\psi(k)=(0,\pi(k)), k \in K$, verify that  $pr\  \phi = \pi = pr\  \psi$, so $0 \to (M, \alpha_M) \stackrel{i} \to (K,\alpha_K) \stackrel{\pi} \to (L, \alpha_L) \to 0$ cannot be a universal central extension.

Lemma \ref{lema 1} ends the proof. \rdg
\bigskip

Classical categories as groups, Lie algebras, Leibniz algebras and other similar ones share the following property: the composition of two central extensions is a central extension, which is absolutely necessary in order to obtain characterizations of the universal central extensions. Unfortunately this property doesn't remain for  the category of Hom-Lie  as the following counterexample \ref{contraejemplo} shows. This problem lead us to introduce the notion of $\alpha$-central extensions in Definition \ref{alfacentral}, whose properties relative to the composition are given in   Lemma \ref{lema 4}.

\begin{Ex} \label{contraejemplo}
Consider the four-dimensional Hom-Lie algebra $(L,\alpha_L)$ with basis $\{a_1,a_2, a_3, a_4\}$, bracket operation given by

$$\left\{\begin{array}{ll} [a_1, a_3] = - [a_3, a_1] = a_4, & [a_1, a_4] = - [a_4, a_1] = a_3,\\

[a_2, a_3] = - [a_3, a_2] = a_1, & [a_2, a_4] = - [a_4, a_2] = a_2. \end{array} \right.$$ (the non-written brackets are equal to zero) and endomorphism $\alpha_L =0$.

Let  $(K,\alpha_K)$ be the five-dimensional Hom-Lie algebra with basis $\{b_1, b_2, b_3, b_4, b_5 \}$, bracket operation given by  $$\left\{\begin{array}{ll} [b_2,b_3] = - [b_3, b_2] = b_1, & [b_2, b_4] = - [b_4,b_2]= b_5,\\

 [b_2,b_5] = - [b_5, b_2] = b_4, & [b_3, b_4] = - [b_4, b_3] = b_2,\\

  [b_3, b_5] = - [b_5, b_3] = b _3. & \end{array} \right.$$ (the non-written brackets are equal to zero) and  endomorphism $\alpha_K =0$.

Obviously $(K,\alpha_K)$ is a perfect Hom-Lie algebra since $K = [K,K]$. On the other hand,  $Z(K,\alpha_K) = <\{b_1 \} >$.

The linear map $\pi : (K, \alpha_K) \to (L,\alpha_L)$ given by $\pi(b_1)=0, \pi(b_2)=a_1,\pi(b_3)=a_2,\pi(b_4)=a_3,\pi(b_5)=a_4$, is a central extension since  $\pi$ is a surjective  homomorphism of  Hom-Lie algebras
and Ker $(\pi) =  <\{b_1 \} > \subseteq Z(K,\alpha_K)$.

Now let us consider the six-dimensional Hom-Lie algebra $(F, \alpha_F)$ with basis $\{e_1, e_2, e_3, e_4, e_5, e_6\}$, bracket operation given by $$\left\{\begin{array}{ll}[e_2, e_3] = - [e_3, e_2] =e_1, & [e_2, e_4] = - [e_4, e_2] =e_1, \\

[e_2, e_5] = - [e_5, e_2] =e_1, & [e_3, e_4] = - [e_4, e_3] = e_2, \\

 [e_3, e_5] = - [e_5, e_3] = e_6, & [e_3, e_6] = - [e_6, e_3] = e_5, \\

 [e_4, e_5] = - [e_5, e_4] = e_3, & [e_4, e_6] = - [e_6, e_4] = e_4, \\

 [e_5, e_6] = - [e_6, e_5] = e_1. & \end{array} \right.$$ (the non-written brackets are equal to zero) and endomorphism $\alpha_F =0$.

The linear map $\rho : (F, \alpha_F) \to (K, \alpha_K)$ given by $\rho(e_1)=0, \rho(e_2)=b_1, \rho(e_3)=b_2, \rho(e_4)=b_3, \rho(e_5)=b_4, \rho(e_6)=b_5$, is a  central extension since  $\rho$ is a surjective homomorphism of Hom-Lie algebras and  Ker $(\rho) =  <\{e_1 \} > = Z(F,\alpha_F)$.

The composition $\pi  \rho : (F, \alpha_F) \to (L, \alpha_L)$ is given by $\pi \rho(e_1) = \pi(0) = 0,  \pi \rho(e_2) = \pi(b_1) = 0, \pi \rho(e_3) = \pi(b_2) = a_1, \pi \rho(e_4) = \pi(b_3) = a_2, \pi \rho(e_5) = \pi(b_4) = a_3, \pi \rho(e_6) = \pi(b_5) = a_4$. Consequently, $\pi \rho : (F, \alpha_F) \to (L, \alpha_L)$ is a surjective homomorphism, but is not a central extension, since $Z(F,\alpha_F) = < \{e_1 \} >$ and Ker $(\pi \rho) =$ $ <\{ e_1, e_2 \} >$, i. e.  Ker $(\pi \rho) \nsubseteq Z(F, \alpha_F)$.
\end{Ex}

\begin{Le} \label{lema 4}
Let $0 \to (M, \alpha_M) \stackrel{i} \to (K,\alpha_K) \stackrel{\pi} \to (L, \alpha_L) \to 0$  and  $0 \to (N, \alpha_N) \stackrel{j} \to (F,\alpha_F) \stackrel{\rho} \to (K, \alpha_K) \to 0$ be  central extensions with   $(K, \alpha_K)$ a perfect Hom-Lie algebra. Then the composition extension $0 \to (P, \alpha_P) = {\rm Ker}\ (\pi \rho)   \to (F,\alpha_F) \stackrel{\pi \rho} \to (L, \alpha_L) \to 0$ is an $\alpha$-central extension.

Moreover, if $0 \to (M, \alpha_M) \stackrel{i} \to (K,\alpha_K) \stackrel{\pi} \to (L, \alpha_L) \to 0$ is a universal  $\alpha$-central extension, then $0 \to (N, \alpha_N) \stackrel{j} \to (F,\alpha_F) \stackrel{\rho} \to (K, \alpha_K) \to 0$ is split.
\end{Le}
{\it Proof.} We must prove that  $[\alpha_P(P), F] = 0$.

Since $(K,\alpha_K)$ is a perfect  Hom-Lie algebra, then every element $f \in F$ can be written as  $f = \displaystyle \sum_i \lambda_i [f_{i_1},f_{i_2}]  + n, n \in N, f_{i_j} \in F, j=1,2$. So, for all $p \in P, f \in F$ we have that
$$[\alpha_P(p), f] = \displaystyle \sum_i \lambda_i \left( [[p,f_{i_1}],\alpha_F(f_{i_2})] + [[f_{i_2},p],\alpha_F(f_{i_1})] \right) + [\alpha_P(p),n] = 0$$
since $[p,f_{i_j}] \in {\rm Ker}\ (\rho) \subseteq Z(F)$.

For the second statement, if $0 \to (M, \alpha_M) \stackrel{i} \to (K,\alpha_K) \stackrel{\pi} \to (L, \alpha_L) \to 0$ is a universal $\alpha$-central extension, then by the first statement, $0 \to (P, \alpha_P) = {\rm Ker}\ (\pi \rho)   \to (F,\alpha_F) \stackrel{\pi \rho} \to (L, \alpha_L) \to 0$ is an  $\alpha$-central extension, then there exists a unique  homomorphism of Hom-Lie algebras  $\sigma : (K,\alpha_K) \to (F,\alpha_F)$ such that $\pi \rho \sigma = \pi$. On the other hand, $\pi \rho \sigma = \pi = \pi Id$ and $(K,\alpha_K)$ is perfect, then Lema \ref{lema 2} implies that $\rho \sigma = Id$.
\rdg

\begin{Th}\label{teorema}\
\begin{enumerate}
\item[a)]  If a central extension  $0 \to (M, \alpha_M) \stackrel{i} \to (K,\alpha_K) \stackrel{\pi} \to (L, \alpha_L) \to 0$ is a universal $\alpha$-central extension, then  $(K,\alpha_K)$ is a perfect   Hom-Lie algebra and every central extension  of $(K,\alpha_K)$ is split.

    \item[b)] Let $0 \to (M, \alpha_M) \stackrel{i} \to (K,\alpha_K) \stackrel{\pi} \to (L, \alpha_L) \to 0$ be a central extension.

     If $(K,\alpha_K)$ is a perfect  Hom-Lie algebra and every  central extension of $(K,\alpha_K)$ is split, then $0 \to (M, \alpha_M) \stackrel{i} \to (K,\alpha_K) \stackrel{\pi} \to (L, \alpha_L) \to 0$ is a universal central extension.

\item[c)] A Hom-Lie algebra $(L, \alpha_L)$ admits a universal central extension if and only if $(L, \alpha_L)$ is perfect.

\item[d)] The kernel of the universal central extension is canonically isomorphic to $H_2^{\alpha}(L)$.

\end{enumerate}
\end{Th}
{\it Proof.}\

\noindent {\it a)} If  $0 \to (M, \alpha_M) \stackrel{i} \to (K,\alpha_K) \stackrel{\pi} \to (L, \alpha_L) \to 0$ is a universal $\alpha$-central extension, then it is a universal central extension by  Remark \ref{rem}, so $(K,\alpha_K)$ is a perfect Hom-Lie algebra by Lemma \ref{lema 3} and every  central extension of  $(K,\alpha_K)$ is split by Lemma \ref{lema 4}.

\bigskip

\noindent {\it b)} Consider a central extension $0 \to (N, \alpha_N) \stackrel{j} \to (A,\alpha_A) \stackrel{\tau} \to (L, \alpha_L) \to 0$. Construct the pull-back extension   $0 \to (N, \alpha_N) \stackrel{\chi} \to (P,\alpha_P) \stackrel{\overline{\tau}} \to (K, \alpha_K) \to 0$, where $P=\{(a,k) \in A \times K \mid \tau(a)=\pi(k) \}$ and  $\alpha_P(a,k)=(\alpha_A(a),\alpha_K(k))$, which is central, consequently is split, i.e. there exists a homomorphism $\sigma : (K,\alpha_K) \to (P,\alpha_P)$ such that $\overline{\tau} \sigma = Id$.

Then $\overline{\pi}  \sigma$, where $\overline{\pi} :  (P,\alpha_P)\to (A,\alpha_A)$ is induced by the pull-back construction, satisfies $\tau \overline{ \pi}  \sigma = \pi$.  Lemma \ref{lema 3} ends the proof.

\bigskip

\noindent {\it c) and d)}   For a Hom-Lie algebra $(L, \alpha_L)$   consider the homology chain complex $C_{\star}^{\alpha}(L)$, which is  $C_{\star}^{\alpha}(L, \mathbb{K})$ where $\mathbb{K}$ is endowed with the trivial  Hom-L-module structure.

As $\mathbb{K}$-vector spaces, let $I_L$ be the  subspace of $L \wedge L$ spanned by the elements of the form $-[x_1,x_2] \wedge \alpha_L(x_3) + [x_1,x_3] \wedge \alpha_L(x_2) - [x_2,x_3] \wedge \alpha_L(x_1), x_1, x_2, x_3 \in L$. That is, $I_L = {\rm Im}\ \left( d_3 : C_3^{\alpha}(L) \to C_2^{\alpha}(L) \right)$.

Now we denote the quotient  $\mathbb{K}$-vector space   $\frac{L \wedge L}{I_L}$  by $\frak{uce}(L)$. Every class $x_1 \wedge x_2 + I_L$ is denoted by $\{x_1,x_2\}$, for all $x_1, x_2 \in L$.

By construction, the following identity holds:
\begin{equation}
\{[x_1,x_2],\alpha_L(x_3)\} + \{[x_2,x_3], \alpha_L(x_1)\} +  \{[x_3,x_1],\alpha_L(x_2)\}=0
\end{equation}
for all $x_1, x_2, x_3 \in L$.

Now $d_2(I_L)=0$, so it induces a $\mathbb{K}$-linear map $u_L : \frak{uce}(L) \to L$, given by $u_L(\{x_1,x_2\})=[x_1,x_2]$. Moreover $(\frak{uce}(L), \widetilde{\alpha})$, where $\widetilde{\alpha} : \frak{uce}(L) \to \frak{uce}(L)$ is defined by $\widetilde{\alpha}(\{x_1,x_2\}) = \{\alpha_L(x_1), \alpha_L(x_2) \}$,  is a  Hom-Lie algebra with respect to the bracket $[\{x_1,x_2\},\{y_1,y_2\}]= \{[x_1,x_2],[y_1,y_2]\}$ and $u_L : (\frak{uce}(L), \widetilde{\alpha}) \to (L,\alpha_L)$ is a homomorphism of Hom-Lie algebras. Actually, Im $u_L = [L,L]$, but $(L,\alpha_L)$ is a perfect Hom-Lie algebra, so $u_L$ is a surjective homomorphism.

From the construction, it follows that  Ker $u_L = H_2^{\alpha}(L)$, so we have the extension
$$0 \to (H_2^{\alpha}(L), \widetilde{\alpha}_{\mid}) \to (\frak{uce}(L), \widetilde{\alpha}) \stackrel{u_L}\to (L,\alpha_L) \to 0$$
which is central, since $[{\rm Ker}\ u_L, \frak{uce}(L)] = 0$, and universal, since for any central extension $0 \to (M,\alpha_M) \to (K,\alpha_K) \stackrel{\pi} \to (L,\alpha_L) \to 0$  there exists the homomorphism of Hom-Lie algebras  $\varphi : (\frak{uce}(L), \widetilde{\alpha}) \to (K,\alpha_K)$ given by $\varphi(\{x_1,x_2\})=[k_1,k_2], \pi(k_i)=x_i, i = 1, 2$, such that $\pi \varphi = u_L$. Moreover, $(\frak{uce}(L), \widetilde{\alpha})$ is a perfect Hom-Lie algebra, so by Lemma \ref{lema 2}, $\varphi$ is unique.
\rdg

\begin{Co} \
\begin{enumerate}

 \item[a)] Let $0 \to (M, \alpha_M) \stackrel{i} \to (K,\alpha_K) \stackrel{\pi} \to (L, \alpha_L) \to 0$ be a universal   $\alpha$-central extension, then $H_1^{\alpha}(K) = H_2^{\alpha}(K) = 0$.

   \item[b)]   Let $0 \to (M, \alpha_M) \stackrel{i} \to (K,\alpha_K) \stackrel{\pi} \to (L, \alpha_L) \to 0$  be a central extension such that $H_1^{\alpha}(K) = H_2^{\alpha}(K) = 0$, then $0 \to (M, \alpha_M) \stackrel{i} \to (K,\alpha_K) \stackrel{\pi} \to (L, \alpha_L) \to 0$  is a universal  central extension.
     \end{enumerate}
\end{Co}
{\it Proof.}

{\it a)} If $0 \to (M, \alpha_M) \stackrel{i} \to (K,\alpha_K) \stackrel{\pi} \to (L, \alpha_L) \to 0$ is a universal  $\alpha$-central extension, then $(K,\alpha_K)$ is perfect by Remark \ref{rem} and Lemma \ref{lema 3}, so  $H_1^{\alpha}(K) = 0$. By Lemma \ref{lema 4} and Theorem \ref{teorema} ( c), d) ) the universal central extension corresponding to $(K,\alpha_K)$ is split, so $H_2^{\alpha}(K) = 0$.

{\it b)} $H_1^{\alpha}(K) = 0$ implies that $(K,\alpha_K)$ is a perfect Hom-Lie algebra.

$H_2^{\alpha}(K) = 0$ implies taht $(\frak{uce}(K),\widetilde{\alpha}) \stackrel{\sim} \to (K,\alpha_K)$.  Theorem \ref{teorema} ( b) ) ends the proof. \rdg

\begin{De} \label{central}
An $\alpha$-central extension $0 \to (M, \alpha_M) \stackrel{i} \to (K,\alpha_K) \stackrel{\pi} \to (L, \alpha_L) \to 0$ is said to be universal if for every central extension $0 \to (R, \alpha_R) \stackrel{j} \to (A,\alpha_A) \stackrel{\tau} \to (L, \alpha_L) \to 0$  there exists a unique homomorphism $\varphi : (K,\alpha_K) \to (A, \alpha_A)$ such that $\tau \varphi = \pi$.
\end{De}

\begin{Pro}
Let $0 \to (M, \alpha_M) \stackrel{i} \to (K,\alpha_K) \stackrel{\pi} \to (L, \alpha_L) \to 0$  and  $0 \to (N, \alpha_N) \stackrel{j} \to (F,\alpha_F) \stackrel{\rho} \to (K, \alpha_K) \to 0$  be central extensions.  If $0 \to (N, \alpha_N) \stackrel{j} \to (F,\alpha_F) \stackrel{\rho} \to (K, \alpha_K) \to 0$ is a universal central extension, then $0 \to (P, \alpha_P) = {\rm Ker} (\pi \rho) \stackrel{\chi} \to (F,\alpha_F) \stackrel{\pi \rho} \to (L, \alpha_L) \to 0$ is an  $\alpha$-central extension which is universal in the sense of Definition \ref{central}.
\end{Pro}
{\it Proof.}  If $0 \to (N, \alpha_N) \stackrel{j} \to (F,\alpha_F) \stackrel{\rho} \to (K, \alpha_K) \to 0$ is a universal central extension, then $(F,\alpha_F)$ and $(K, \alpha_K)$ are perfect Hom-Lie algebras by Lemma \ref{lema 3}.

On the other hand, $0 \to (P, \alpha_P) = {\rm Ker} (\pi  \rho) \stackrel{\chi} \to (F,\alpha_F) \stackrel{\pi  \rho} \to (L, \alpha_L) \to 0$ is an $\alpha$-central extension by Lemma \ref{lema 4}.

In order to obtain the universality, for any central extension $0 \to (R, \alpha_R) \to (A,\alpha_A) \stackrel{\tau} \to (L, \alpha_L) \to 0$  construct the pull-back extension corresponding to $\tau$ and $\pi$, $0 \to (R,\alpha_R) \to (K \times_L A, \alpha_K \times \alpha_A) \stackrel{\overline{\tau}}\to (K, \alpha_K) \to 0$. Since $0 \to (N, \alpha_N) \stackrel{j} \to (F,\alpha_F) \stackrel{\rho} \to (K, \alpha_K) \to 0$ is a universal central extension, then there exista a unique homomorphism $\varphi : (F, \alpha_F) \to (K \times_L A, \alpha_K \times \alpha_A)$ such that $\overline{\tau} \cdot \varphi = \rho$. Then the homomorphism $\overline{\pi} \cdot \varphi$ satisfies that $\tau \cdot \overline{\pi} \cdot \varphi = \pi \cdot \rho$ and it is unique by Lemma \ref{lema 2}.

\rdg

\bigskip

\centerline{\bf Acknowledgements}
First and second  authors  were supported by  Ministerio de
Educaci\'on y Ciencia (Spain), Grant MTM2009-14464-C02 (European
FEDER support included) and by Xunta de
Galicia, Grant Incite09 207 215 PR.

\bigskip

\begin{center}

\end{center}

\end{document}